\documentclass[12pt]{article}       
\usepackage[english]{babel}
\usepackage{amsmath,amssymb,amsfonts,makeidx,amssymb}
\usepackage{graphicx}
\usepackage{lscape}
\usepackage{todonotes}
\usepackage{natbib}
\usepackage{color}%
\usepackage[normalem]{ulem}
\newcommand{\ind}{\mathbb{I}}

\newtheorem{lem}{Lemma}

\newtheorem{theo}{Theorem}
\newtheorem{cor}{Corollary}

\newtheorem{prop}{Proposition}
\newtheorem{rem}{Remark}

\begin{document}

	\begin{center}
		\Large \bf A combined strategy for multivariate density estimation\\
	\end{center}
	\normalsize
	
	\
	
	\begin{center}
		Alejandro Cholaquidis$^a$, Ricardo Fraiman$^a$, Badih Ghattas$^b$ \\ and Juan Kalemkerian$^a$\\
		$^a$ Universidad de la Rep\'ublica, Uruguay\\
		$^b$ Aix Marseille Universit\'e, CNRS, Marseille, France.\\
	\end{center}

\begin{abstract}
Non-linear aggregation strategies have recently been proposed in response to the problem of how to 
combine, in a non-linear way, estimators of the regression function (see for instance \cite{biau:16}), 
classification rules (see \cite{ch:16}), among others. Although there are several linear strategies to aggregate density 
estimators, most of them are hard to compute (even in moderate dimensions). 
Our approach aims to overcome this problem by estimating the density at a point $x$ using not just 
sample points close to $x$ but in a neighborhood of the (estimated) level set $f(x)$. We show, both 
theoretically and through a simulation study, that the mean squared error of our proposal is smaller than
that of the aggregated densities. A Central Limit Theorem is also proven.

\end{abstract}

\section{Introduction}
\begin{figure}[h]
	\begin{center}
		\includegraphics[scale=.3]{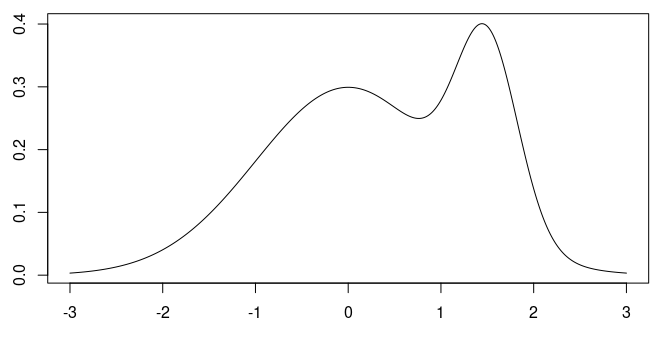} 
				\includegraphics[scale=.3]{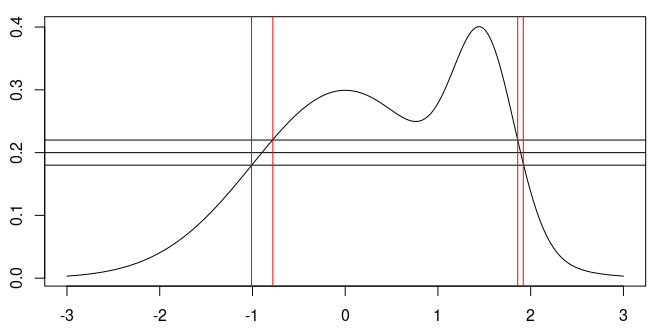} 
		\caption{Left: a density whose concentration mass varies significantly within its support. Right: the $0.2$-neighborhood for the level $f(x)=0.2$ is given by the union of the intervals $I_1=[-1.01,-0.78]$ and $I_2=[1.86,1.92]$.}
		\label{ex1}
	\end{center}
\end{figure}

\begin{figure}[t!]
	\begin{center}
		\includegraphics[scale=.4]{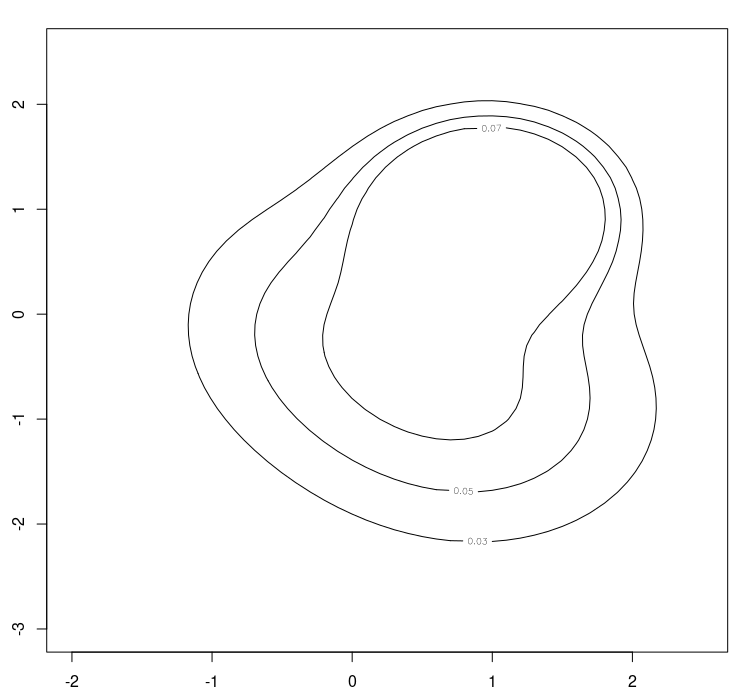}
		\caption{$0.02$-neighbourhood of a level set for a mixture of three bi-variate Gaussian distributions, the first with mean $(0,0)$ and covariance matrix $\Sigma=diag(1,1)$, the second with mean $(1,1)$ and covariance matrix $\Sigma=diag(1/2,1/2)$, the last with mean $(1,-1)$ and covariance matrix $\Sigma=diag(1,1)$, $f(x)=0.05$}
		\label{ex2}
	\end{center}
\end{figure}

\indent Density estimation is still an important and active area of research that has many 
statistical applications, particularly in supervised and unsupervised learning, 
see for instance the recent book by \cite{cd:18}. Although  this is a well-studied subject, 
when the data belongs to high or even moderate dimensions, such as
$\mathbb{R}^2$ or $\mathbb{R}^4 $, this becomes a difficult problem due to the well-known curse
of dimensionality. This is also the case for non-parametric regression. 
For this last problem, \cite{biau:16} propose a non-linear aggregation method that is very close in spirit to our approach. In \cite{ch:16}, the authors propose a similar idea for classification. To tackle this problem, we introduce a new non-linear aggregation method that is well designed for moderate dimensions. Our approach is based on two main ideas:
\begin{itemize}
	\item[1)] The first idea is to compute the estimator of $f$ at the point $x$ 
	using an estimator of a $\epsilon$-neighborhood of the level set, i.e,  
	$ \{y:\vert f(y)-f(x)\vert  \leq \epsilon \}\equiv B^*( \epsilon,x),$
	instead of a neighborhood of the point $x$, see the right-hand panel of Figure \ref{ex1}  and also see Figure \ref{ex2}. Roughly speaking, under the unrealistic case where $B^*(\epsilon,x)$ is known, the estimator that we  propose behaves as if the data were in one dimension. In general $B^*(\epsilon,x)$ is unknown, consequently a loss of efficiency will appear, which is related to the estimation of the $\epsilon$-neighborhood. 

	\item[2)] The second idea is to perform a nonlinear aggregation method to combine several estimators. This  will improve the behavior when, for instance, the underlying true density $f$ is not unimodal and the concentration of mass varies significantly within its support, see Figure \ref{ex1}.

\end{itemize}

Similar ideas have previously been considered for density estimation and non-parametric regression. With respect to 1), a related approach can be found in \cite{flm:97}, where it is assumed that the density has a particular shape given by the composition of a univariate density with a depth. The particular case of ellipsoidal density has also been considered in \cite{stute}.  In our setup, no particular structure is required to the multivariate density. 

 Starting from the seminal work by \cite{brei:96}, many linear aggregation methods have been developed, see for instance \cite{le:06,rt:07,bg:13,be:17} and the references therein.

The rest of this paper is organised as follows. In Section \ref{notat}, we introduce the notation 
and main definitions used through the manuscript. In Section \ref{comb}, we define a 
nonlinear aggregation estimator $\hat{f}_{\text{agg}}$ that is based on a family of density 
estimators $f_1, \ldots, f_M$, which requires perfect match of the sets 
$ \{y:\vert f_j(y)-f_j(x)\vert  \leq \epsilon \}$ for $j=1,\ldots,M$, 
see \eqref{beeps}.  It  can also be relaxed to a partial matching because it will be  
defined later.  In Subsection \ref{opti}, we prove that the aggregated strategy 
is asymptotically optimal in the sense that it behaves as well as the best density 
estimator within the family. In Subsection \ref{cons} we prove consistency 
in $L^2$, under mild regularity conditions on $f$. A Central Limit Theorem 
is proven in Subsection \ref{central} for the case $M=1$. Lastly, 
in Section \ref{simu} we perform some simulations in dimensions $2$ and $4$, which illustrate the good performance of our approach.

\section{Notation} \label{notat}

Let us consider $\mathbb{R}^d$ endowed with the $d$-dimensional Lebesgue measure $\mu$.
For $r>0$, $\mathcal{B}(x,r)$ denotes the open ball of radii $r>0$, 
and $\omega_d=\mu(\mathcal{B}(0,1))$. Given $A\subset \mathbb{R}^d$, we will denote
by $B(A,\epsilon)$ the parallel set of radius $\epsilon$ of $A$, that is 
$B(A,\epsilon)=\{y\in \mathbb{R}^d: d(y,A)<\epsilon\}$ where 
$d(y,A)=\inf_{a\in A} \|y-a\|$, and $\|\cdot\|$ denotes the Euclidean norm. Given a 
kernel function $K:\mathbb{R}^d\rightarrow \mathbb{R}^+$ we say that $K$ is regular
if there exists $0<c_1<c_2<\infty$ such that $c_1\mathbb{I}_{\mathcal{B}(0,1)}(x)\leq
K(x)\leq c_2\mathbb{I}_{\mathcal{B}(0,1)}(x)$, where $\mathbb{I}_{A}$ 
stands for the indicator function of the set $A$. 
We will denote $K_h(x)=K(x/h)$ and  $B^*(\epsilon,x)=\{y:|f(x)-f(y)|<\epsilon\}$.

\section{The combined estimator} \label{comb}

Throughout this manuscript, we will assume that $f$ is a density, bounded from above, 
such that  $f(X)\in L^2$. Let $\mathcal{D}_{n}=\{X_1,\dots,X_n\}$ be iid 
random vectors with the same distribution $f$ as $X$. We split $\mathcal{D}_n$ into two 
disjoint subsets, namely $\mathcal{D}_{k}=\{X_1,\dots,X_k\}$ and $\mathcal{E}_l=
\{X_{k+1},\dots,X_{n}\}$ with $l=n-k$.
Let $\mathbf{f_k}(x)= (f_1(x),\dots,f_M(x))$ be $M$ density estimators computed 
with the first sample $\mathcal{D}_{k}$. 

For $\epsilon>0$, we define the combined neighborhood of radius $\epsilon$, $B(\epsilon,x)$, of a given point $x$ to be 
\begin{equation} \label{beeps}
B(\epsilon,x)=  \left\{y\in \mathbb{R}^d:  \bigcap_{m=1}^M |f_m(y) - f_m(x)| < \epsilon \right\}.
\end{equation}
Let us consider the estimator of $P_X(B(\epsilon,x))$, given by
\begin{align}\label{n1}
N(\epsilon,x) =& \# \left\{\bigcap_{m=1}^M\left\{X_i \in \mathcal{E}_l:   |f_m(X_i) - f_m(x)| < \epsilon \right\}\right\} /l\\
 =&\frac{1}{l} \sum_{j=1}^l \prod_{m=1}^M \mathbb{I}_{[0,\epsilon)} \big(|f_m(X_{j+k})-f_m(x)|\big).
\end{align}
Lastly, the aggregated density estimator is defined as 
\begin{equation}\label{eq1}
\hat{f}_{\text{agg}}(x)= \frac{N(\epsilon,x)}{\mu(B(\epsilon,x))}= \frac{\sum_{j=1}^l \mathbb{I}_{B(\epsilon,x)}(X_{j+k})}{l\mu(B(\epsilon,x))}.
\end{equation}

\subsection{A smoothed approach}

Instead of the indicator function used in \eqref{n1}, we can use a one dimensional kernel $K$. Define,

\begin{equation*}
\tilde{N}(\epsilon,x) = \frac{1}{l}\sum_{j=1}^l \prod_{m=1}^M K\Big[\frac{f_m(X_{k+j})-f_m(x)}{\epsilon}\Big],
\end{equation*}
where  $K$ fulfils, 
\begin{equation}\label{eq00}
\lim _{\epsilon \to 0} \frac{1}{\mu(B^*(\epsilon,X))}\int_{B^*(\epsilon,x)} K\Big[\frac{f(t)-f(X)}{\epsilon}\Big]^Mdt = 1 \quad a.s.
\end{equation}
Then the alternative aggregated estimator is defined as, 
\begin{equation}\label{eqest2}
\tilde{f}_{\text{agg}}(x)= \frac{\tilde{N}(\epsilon,x)}{\mu(B(\epsilon,x))}.
\end{equation}

\subsection{An alternative approach}

Let $\epsilon>0$ and $0\leq \eta<1$, define the $\eta$-neighborhood of radius $\epsilon$, $B^\eta(\epsilon,x)$, of a given point $x$ to be 
$$B^\eta(\epsilon,x)=  \left\{y\in \mathbb{R}^d:  \frac{1}{M}\sum_{m=1}^M \mathbb{I}_{\{|f_m(y) - f_m(x)| < \epsilon\}}\geq 1-\eta \right\}.$$
Observe that the for $\eta=0$ we get $B^\eta(\epsilon,x)=B(\epsilon,x)$. We define the $\eta$-density estimator,  $\hat{f}_{\text{agg},\eta}(x)$ as in \eqref{eq1} replacing $B(\epsilon,x)$ with $B^\eta(\epsilon,x)$. Regarding \eqref{eqest2} we can define
\begin{equation}\label{tildefeta}
\tilde{f}_{\text{agg},\eta}(x)= \frac{1}{l}\sum_{j=1}^l \mathbb{I}_{\{X_{k+j}\in B^\eta(\epsilon,x)\}} \prod_{m=1}^M  K\Big[\frac{f_m(X_{k+j})-f_m(x)}{\epsilon}\Big].
\end{equation}

\subsection{Optimality} \label{opti}
The following proposition (which is the analogous for our setup of  Proposition 2.1 in \cite{biau:16}), states that the combined 
estimator behaves as well as the best density estimator, except for the second term, which will be proven to converge to 0 (see Theorem \ref{thcons}).

\begin{prop} \label{opt} With the notation introduced previously 
	\begin{equation} \label{prop1}
	\mathbb{E}|\hat{f}_{\emph{agg}}(X)-f(X)|^2\leq \min_{m=1,\dots,M} \mathbb{E}|f_m(X)-f(X)|^2+\mathbb{E}|\hat{f}_{\emph{agg}}(X)-T(\mathbf{f_k}(X))|^2,
	\end{equation}
where $X$ is independent of $\mathcal{D}_n$ and $T(\mathbf{f_k}(X))=\mathbb{E}(f(X)|\mathbf{f_k}(X)).$
\end{prop}

\begin{rem} It is easy to see that Proposition \ref{opt} holds for $\tilde{f}_{\text{agg}}(X)$.
\end{rem}

Now we will state two Lemmas (whose proofs are given in the Appendix), the first  proves that the theoretical estimator $T(\mathbf{f_k}(X))$ converges in
 $L^2$ to $f(X)$, as $k\rightarrow \infty$. The second  proves that under point-wise consistency of the density estimators,  for all $\epsilon>0$, with probability one $\mu(B(\epsilon,x))\rightarrow \mu(B^*(\epsilon,x))$ as $k\rightarrow \infty$ for almost all $x$. To do this, let us introduce the following condition,
 
\begin{itemize}
	\item[K1] A random variable $X$ with distribution $P_X$ and density $f$ fulfils $K1$, if  $\mathbb{P}(f(X)=a)=0$ for all $a\in \mathbb{R}$.
\end{itemize}

\begin{lem} \label{L2conv} Under $K1$, if $\{f_i\}_i$ is any sequence of functions (possibly random) such $\lim_{i\rightarrow \infty}f_i(X) = f(X)$ a.s  then,
$$\lim_{i\rightarrow \infty} \mathbb{E}\big|\mathbb{E}[f(X)|f_i(X)]-f(X)\big|^2=0.$$
\end{lem}

\begin{lem}\label{lem} Let $X$ be random variable with distribution $P_X$ whose density $f$ is continuous. Let $\mathcal{D}_k=\{X_1,\dots,X_k\}$ be an iid sample of $X$ and $f_1,\dots,f_M$ be continuous density estimators (built from $\mathcal{D}_k$),  such that for all $m=1,\ldots,M$, $|f_m(x)-f(x)|\rightarrow 0$ a.s., as $k\rightarrow \infty$ for almost all $x$ w.r.t to $\mu$.  Let $\epsilon>0$, then for all $x$ such that
	\begin{itemize}
		\item $f_m(x)\rightarrow f(x)$ for $m=1,\ldots,M$, a.s., as $k\rightarrow \infty$.
		\item $\mu[B^*(\epsilon+\gamma,x)\setminus B^*(\epsilon-\gamma,x)]\rightarrow 0$ as $\gamma\rightarrow 0$.
		\item $\overline{B^*(\epsilon,x)}$ is compact, and $\overline{B(\epsilon,x)}$ is compact a.s.
		\end{itemize}  we have 
	\begin{equation} \label{lem1}
	\mu(B(\epsilon,x))\rightarrow \mu(B^*(\epsilon,x))\quad a.s., \text{ as } k\rightarrow \infty,
	\end{equation}
	and
	\begin{equation}\label{lem2}
	P_X(B(\epsilon,x))\rightarrow P_X(B^*(\epsilon,x))\quad a.s., \text{ as }k\rightarrow \infty.
	\end{equation}
\end{lem}

\subsection{Consistency} \label{cons}

Because the first term  in the right-hand side of \eqref{prop1} does not depend on $l$ and converges to $0$ if at least one of the density estimators  is mean square error consistent, to prove the consistency (taking limit first in $l$ and second in $\epsilon$) for the aggregated estimator, we only need to prove that the second term in the right-hand side of \eqref{prop1} converges to $0$ in mean square error. This is done in the following Theorem, under mild regularity restrictions on $P_X$, as well as point-wise convergence for the density estimators  and uniform equicontinuity. Recall that a sequence of functions $\{g_k\}_k$ is said to be uniformly equicontinuous if for all $\epsilon>0$ there exists $\delta=\delta(\epsilon)$ such that for all $k$, $|g_k(x)-g_k(y)|<\epsilon$, whenever $\|x-y\|<\delta$. All of the proofs of this section are given in the Appendix.

\subsubsection{Assumptions} 

We will consider the following set of assumptions

\begin{itemize}
	\item[H1] The density estimators $f_1,\dots,f_M$ based on a sample $\mathcal{D}_k$ fulfils H1 if with probability one, the sequences $\{f_1\}_k,\ldots,\{f_M\}_k$  are uniformly equicontinuous and the $\delta=\delta(\epsilon)$ 	of the uniform equicontinuity is bounded from below by $\delta_0(\epsilon)>0$. 
	\item[H2] The density estimators $f_1,\dots,f_M$ based on a sample $\mathcal{D}_k$ fulfils H2 if for almost all $x$ w.r.t. $\mu$, $f_j(x)\rightarrow f(x)$, a.s., for all $j=1,\dots,M$ as $k\rightarrow \infty$.
	 
\end{itemize}

\begin{theo}\label{thcons}  
Let us assume K1, H1 and H2. We assume also that, for all $x$ such that $f_m(x)\rightarrow f(x)$ for all $m=1,\ldots,M$, there exists $\epsilon_0(x)$ such that for all $0<\epsilon<\epsilon_0(x)$, the set $\overline{B^*(\epsilon,x)}$ is compact, the set $\overline{B(\epsilon,x)}$ is compact a.s., and $\mu[B^*(\epsilon+\gamma,x)\setminus B^*(\epsilon-\gamma,x)]\rightarrow 0$ as $\gamma\rightarrow 0$. Let $k=k(l)\rightarrow \infty$ 
as $l\rightarrow \infty$, then, 
\begin{equation*}
\lim_{\epsilon \rightarrow 0} \lim_{l\rightarrow \infty}\mathbb{E}|\hat{f}_{\emph{agg}}(X)-T(\mathbf{f_k}(X))|^2= 0.
\end{equation*}
\end{theo}

\begin{theo}\label{th12} Under the hypotheses of Theorem \ref{thcons}.  If $K$ is a kernel function, bounded from above by $c_2<\infty$, that fulfils \eqref{eq00}, then
	\begin{equation*}
\lim_{\epsilon \rightarrow 0} \lim_{l\rightarrow \infty}\mathbb{E}|\tilde{f}_{\emph{agg}}(X)-T(\mathbf{f_k}(X))|^2= 0.
	\end{equation*}
	
\end{theo}

\begin{rem}
	\begin{itemize} 
		\item[1)] Corollary 1  \cite{eima:05} proves that if $f$ is uniformly continuous (with some regularity conditions  on the
		kernel $K$), then the multidimensional kernel density estimator converges almost surely, uniformly, by choosing a suitable bandwidth. 
		It is easy to see that this entails the required uniform equicontinuity on the estimators. 
		\item[2)] Following the same ideas used to prove Theorem \ref{thcons}, it can be proven that
		$\lim_{\epsilon \rightarrow 0}\lim_{l \rightarrow \infty} \mathbb{E}(\tilde{f}_{\text{agg}}(X)-f(X))^2=0$ (see Theorem \ref{th12} in Appendix).
\end{itemize}
\end{rem}

 If the density $f$ is bounded from below by a positive constant, we have the following direct corollary.
\begin{cor}\label{cor1} Under the hypotheses of Theorem \ref{thcons}, if  in addition the 
density $f$ fulfils  that there exists $C$ and $A$ such that $0<A\leq f(x)\leq C<\infty$ for all $x$, then, 
	\begin{equation*}
 \lim_{\epsilon \rightarrow 0} \lim_{l\rightarrow \infty}\mathbb{E}\int_{\mathbb{R}^d}|\hat{f}_{\emph{agg}}(x)-f(x)|^2dx= 0.
	\end{equation*}
\end{cor}

\subsection{A central limit theorem} \label{central}


The following theorem states that a central limit theorem for $\hat{f}_{\text{agg}}(x)$ holds, when the limit is taken first as $k\rightarrow \infty$ and second as $l\rightarrow \infty$.

\begin{theo} \label{tcl} Let $\epsilon=\epsilon_l\rightarrow 0$ such that $l\epsilon_l^2\rightarrow 0$. Then, for all $x$ such that $f(x)>0$ and 
	\begin{itemize}
		\item $\mu(\{y:f(x)=f(y)\})=0$
		\item Exists $\epsilon_0>0$, such that  $\overline{B^*(\epsilon',x)}$  is compact, and $\overline{B(\epsilon',x)}$ is compact a.s.
		for all $\epsilon'<\epsilon_0$.
		\item  Exists $\epsilon_0>0$, such that $\mu[B^*(\epsilon'+\gamma,x)\setminus B^*(\epsilon'-\gamma,x)]\rightarrow 0$ as
		$\gamma\rightarrow 0$ for all $\epsilon'<\epsilon_0$.
		\item  $\mu(B^*(\epsilon,x))l\rightarrow \infty$ as $l \rightarrow \infty$.
		\item $f_m(x)\rightarrow f(x)$ for all $m=1,\dots,M$ as $k \rightarrow \infty$.
	\end{itemize}
We have,
\begin{equation}\label{clt}
\lim_{l\rightarrow \infty}\lim_{k\rightarrow \infty} \sqrt{\mu(B^*(\epsilon,x))l }\bigg[\hat{f}_{\emph{agg}}(x)-f(x)\bigg]\stackrel{d}{=} N(0,f(x)).
\end{equation}
	\end{theo}

\begin{rem}
The previous theorem depends on the calculus of $\mu(B^*(\epsilon,x))$, which is in general unknown. However, in some cases it can be estimated, by means of a Monte-Carlo method, using a uniformly consistent estimator $f_n$ of $f$ and a sample of uniformly distributed random variables on a box containing the set $B^*(\epsilon,x)$.\\

For the special case of spherical densities (i.e., $f(x)=h(\|x\|^2)$ for some $h:\mathbb{R}\rightarrow \mathbb{R}$), the limit of $\mu(B^*(\epsilon,x))/\epsilon$ can easily be derived, as  is proven in the following proposition. 
\end{rem}

\begin{prop} \label{prop2} Let $f$ be a spherical density such that $h$ is strictly decreasing and $h'$ is continuous on a neighbourhood containing $\|x\|^2$, then, for all $x$ such that $f(x)>0$, and $\|\nabla f(x)\|>0$,
	$$\lim_{l\rightarrow \infty} \frac{\mu(B^*(\epsilon,x))}{2\epsilon}= \frac{2\pi^{d/2}\|x\|^{d-1}}{\Gamma(\frac{d}{2})\|\nabla f(x)\|},$$
	where $\Gamma$ is Euler's gamma function.
\end{prop}

\section{Models used for the simulations} \label{simu}

First, we performed a simulation study  to assess, in terms of the mean square error, the proposed aggregation strategy. Second, we evaluate the departure from normality in Theorem \ref{tcl}. Five different distributions were considered:
\begin{itemize}
	\item[1] Beta, with density
	$\left(\frac{\Gamma(\alpha+\beta)}{\Gamma(\alpha)\Gamma(\beta)}\right)^d(x_1\cdots x_d)^{\alpha-1}(1-x_1)^{\beta-1}\cdots (1-x_d)^{\beta-1}.$
	 \item[2] Normal, with mean $0$ and variance $\Sigma=\text{diag}(\sigma_1^2,\dots,\sigma_d^2)$ is a diagonal matrix.
	 \item[3] Weibull, with density 
	 $\left(\frac{k}{\lambda^k}\right)^d(x_1\cdots x_d)^{d(k-1)}\exp\Big(-\sum^d_{i=1} (x_i/\lambda)^k\Big)$.
	 \item[4] Convex combination of two bi-variate normal distributions with the same covariance matrix $\Sigma$: $(1/2)N(\mu_1,\Sigma)+(1/2)N(\mu_2,\Sigma)$
	 where $\Sigma=\Sigma_1$ given below. 
	 \item[5] Convex combination of two bi-variate normal distributions: $(1/2)N(\mu_1,\Sigma_1)+(1/2)N(\mu_2,\Sigma_2)$
	 where 
	 $$\Sigma_1= \begin{bmatrix}
	 \sigma_1^2 &\rho\\
	 \rho & \sigma_2^2
	 \end{bmatrix}\quad \quad \text{ and } \quad\Sigma_2= \begin{bmatrix}
	 \sigma_2^2 &-\rho\\
	 -\rho & \sigma_1^2
	 \end{bmatrix}.$$
\end{itemize}

To build the estimator $\hat{f}_{\text{agg}}$ we considered five kernel-based density estimators $f_{k,\gamma_1},\dots,f_{k,\gamma_5}$ computed with different bandwidth $\gamma_1,\dots,\gamma_5$. The bandwidths were chosen as follows: first we compute the leave-one-out cross validation bandwidth $hcv$ based on a sample of size $k$. This value is kept fixed along the
replicates. Then, we fix $\gamma_1=0.9\times hcv$, $\gamma_2=0.95\times hcv$, $\gamma_3=hcv$, $\gamma_4=1.05\times hcv$ and $\gamma_5=1.1\times hcv$. We choose $k=l=2000$ for $d=2$ and $k=l=4000$ for $d=4$. Let us denote $hcvu$ the leave-one-out cross validation bandwidth based on the whole sample. The parameter $\epsilon_l$ was selected as follows: first, we compute the five kernel-based density estimators $f_{k+l,\tilde{h}_1},\dots,f_{k+l,\tilde{h}_5}$ based on the whole sample $\mathcal{D}_k\cup \mathcal{E}_l$, with bandwidth $\tilde{h}_1=0.9\times hcvu$, $\tilde{h}_2=0.95\times hcvu$, $\tilde{h}_3=hcvu$, $\tilde{h}_4=1.05\times hcvu$ and $\tilde{h}_5=1.1\times hcvu$,  we then compute the average of them; i.e, $\overline{f}(x)=(f_{k+l,\tilde{h}_1}+\ldots +f_{k+l,\tilde{h}_5})/5.$

Finally, $\epsilon_l$ is the value that minimize $\|\hat{f}_{\text{agg}}-\overline{f}\|_2$ (where $\|\cdot\|_2$ denotes the $L^2$ norm).

\begin{figure}[h]
	\begin{center}
		\includegraphics[scale=.25]{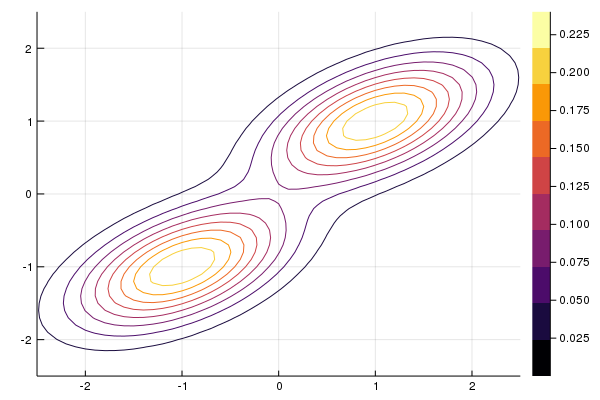}  
		 		\includegraphics[scale=.25]{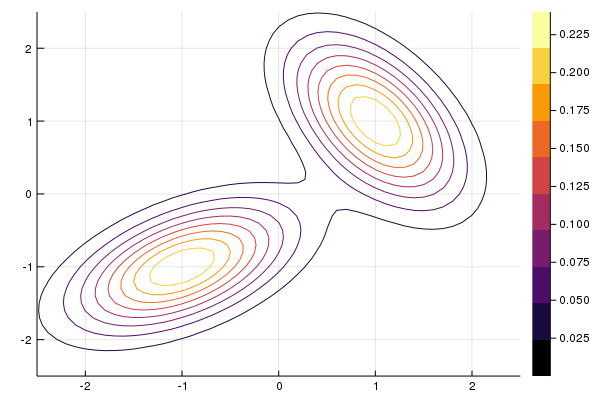}
		\caption{Left-hand panel: 10 level sets for the density of model 4. Right-hand panel: eight level sets for the density of model 5. In both models, we used $\rho=.2$, $\sigma_1=.5$ and $\sigma_2=.3$.}
		\label{conv}
	\end{center}
\end{figure}

The measures $\mu(B(\epsilon_l,x))$ are computed by Monte-Carlo method using $20000$ and $40000$ uniformly distributed random variables in dimensions $2$ and $4$, respectively. Two different kernels where considered: the Epanechnikov kernel (denoted by E), and the Gaussian kernel (denoted by G). The whole procedure is repeated 100 times. We report  $\|\hat{f}_{\text{agg}}-f\|_2$, estimated from a test sample, uniformly distributed over a rectangle in $\mathbb{R}^2$, or $\mathbb{R}^4$. \\
	Figure \ref{conv} shows the level sets for the density of model 4 (left panel) and model 5 (right panel).

\begin{table}[h]
\caption{$L_2$ error over $100$ replicates for model 1. For the test sample, we used $2000$ uniformly distributed points on $[0,1]^2$. The measure of $B(\epsilon,x)$ is estimated using $20000$ uniformly distributed points on $[0,1]^2$ for $d=2$, and $40000$ for $d=4$.}
\label{beta}
	\tiny
	\begin{center}
		\begin{tabular}{|c|c|c|c|c|}
	\hline
	& \multicolumn{2}{c|}{$\alpha=1.5=\beta$} &	\multicolumn{2}{c|}{$\alpha=2.5=\beta$} \\
	\hline
                            &$d=2$            & $d=4$          &$d=2$ &$d=4$\\
     	\hline
	$n,k$     			         & 2000        & 4000          &2000 &4000\\
Kernel                       & G              &G               & G   &   G\\
\hline
${\hat{f}_{\text{agg}}}$    & \textbf{0.090}  & \textbf{0.185} &0.118 			&0.269  \\
$f_{k,\gamma_1}$             & 0.111           &0.240           &0.122 			&0.305 \\
$f_{k,\gamma_2}$            & 0.110           &0.233           &0.123          &0.300\\
$f_{k,\gamma_3}$                 & 0.111            &0.231          &0.125          &0.301 \\
$f_{k,\gamma_4}$            & 0.113           &0.232           &0.128		    &0.307   \\
$f_{k,\gamma_5}$             & 0.115           &0.235           &0.131 			&0.316  \\
$f_{n,\gamma_1}$           & 0.093           &0.200           &0.100		    &0.256   \\
$f_{n,\gamma_2}$          & 0.095            &0.201          &0.103 			&0.261\\
$f_{n,hcvu}$              & 0.092           &0.200           &\textbf{0.098} &\textbf{0.255} \\
$f_{n,\gamma_4}$          & 0.101           &0.211           &0.112          &0.282 \\
$f_{n,\gamma_5}$           & 0.104           &0.218           &0.117          & 0.297   \\
 \hline      	
\end{tabular}
\hfill
\hfill	\begin{tabular}{|c|c|c|c|c|}
	\hline
	& \multicolumn{2}{c|}{$\alpha=1.5=\beta$} &	\multicolumn{2}{c|}{$\alpha=2.5=\beta$} \\
	\hline
	&$d=2$            & $d=4$            &$d=2$ &$d=4$\\
	\hline
	$k,l$     			         & 2000           & 4000       	      &2000 &4000\\
	Kernel                       & E              &  E      		  & E   &   E\\
	\hline
	${\hat{f}_{\text{agg}}}$    & \textbf{0.090}   & 0.141           &\textbf{0.099} &\textbf{0.229} \\
	$f_{k,\gamma_1}$             & 0.107           & 0.158            &0.162          &0.261  \\
	$f_{k,\gamma_2}$            & 0.104           & 0.151            &0.161          &0.252    \\
	$f_{k,\gamma_3}$                 & 0.102           & 0.149            &0.162          &0.250  \\
	$f_{k,\gamma_4}$            & 0.102           & 0.151            &0.164          &0.256         \\
	$f_{k,\gamma_5}$             & 0.105           & 0.158            &0.170          &0.268       \\
	$f_{n,\gamma_1}$           & 0.102            & 0.145           &0.157          &0.234         \\
	$f_{n,\gamma_2}$          & 0.099           & 0.137            &0.156          & 0.229       \\
	$f_{n,hcvu}$              & 0.098           & \textbf{0.136}   &0.157          &0.230        \\
	$f_{n,\gamma_4}$          & 0.099           & 0.142            &0.161          &0.241     \\
	$f_{n,\gamma_5}$           & 0.102           & 0.150            &0.167          &0.256    \\
	\hline      	
	\end{tabular}
	\end{center}
\end{table}
\begin{table}[h]
\caption{$L_2$ error over $100$ replications for model 2. The test sample consist of 2000 uniformly distributed points on $[0,4]^2$ in dimension 2 and 
	4000 uniformly distributed points on $[0,5]^4$ in dimension 4.}
	\label{wei}
	\tiny	
	\begin{center}
\begin{tabular}{|c|c|c|c|}
		\hline
									 &  $\lambda=1,k=1$&  $\lambda=1, k=0.5$&  $\lambda=1, k=1$ \\
		\hline
									 &$d=2$,           &  $d=2$             &  $d=4$ \\
		\hline
		$n,k$     			         & 2000            &  2000              &   40000\\
		Kernel                       & E               &  E                 &   E  \\
		\hline
		${\hat{f}_{\text{agg}}}$     & \textbf{0.069}  & \textbf{0.653}    &  \textbf{0.009}\\
		$f_{k,\gamma_1}$              & 0.119           & 0.678             &   0.026\\
		$f_{k,\gamma_2}$             & 0.113           & 0.677   			&  0.024 \\
		$f_{k,\gamma_3}$                  & 0.108           & 0.676             &   0.021\\
		$f_{k,\gamma_4}$             & 0.103           & 0.675             &   0.019\\
		$f_{k,\gamma_5}$              & 0.098           & 0.674             &   0.018\\
		$f_{n,\gamma_1}$            & 0.086           & 0.672             &   0.018\\
		$f_{n,\gamma_2}$           & 0.082           & 0.671             &   0.017\\
		$f_{n,hcvu}$               & 0.078           & 0.694             &   0.020\\
		$f_{n,\gamma_4}$           & 0.074           & 0.670             &   0.014\\
		$f_{n,\gamma_5}$            & 0.072           & 0.670             &   0.012\\
		\hline      	
	\end{tabular}
	\end{center}
\end{table}
\begin{table}[h]
\caption{$L_2$ error for model 2 over $100$ replicates using Epanechnikov's kernel. In $\mathbb{R}^2$ $\Sigma=diag(\sigma_1^2,\sigma_2^2)$, and  in $\mathbb{R}^4$  $\Sigma=diag(\sigma_1^2,\sigma_2^2,\sigma_3^2,\sigma_4^2)$. The test consist of 2000 uniform on $[-3,3]\times[-2,1.5]$ for the first column and on  $[-2.5,2.5]^2$ for the second column. For dimension $4$, the test sample is uniformly distributed on $[-3,3]\times [-1.5,1.5]^2\times [-3,3]$.}
	\label{normal}
	\tiny
	\begin{center}
	\begin{tabular}{|c|c|c|c|}	
		\hline
	                             	 & d=2                      & d=2                         & $d=4$\\
		\hline
		                             &$\sigma_1=1,\sigma_2=0.25$  & $\sigma_1=1,\sigma_2=0.1$  & $\sigma_1=1=\sigma_4$,\\
		                             &                           &                            & $\sigma_2=.5=\sigma_3$\\
		                 \hline		                          
		$n=k$ 	                     & 2000                      & 2000                       & 4000\\
		Kernel	                     & E                         &  E                         & E\\
		\hline
		${\hat{f}_{\text{agg}}}$     & \textbf{0.023}        	&\textbf{0.082}           & \textbf{0.008}  \\
		$f_{k,\gamma_1}$              & 0.034 			        &  0.126                  & 0.050  \\
		$f_{k,\gamma_2}$             & 0.033                    &  0.120  		          & 0.045   \\ 
		$f_{k,\gamma_3}$                  & 0.032                    &  0.114                  & 0.041  \\
		$f_{k,\gamma_4}$             & 0.030                    &  0.108    		      & 0.037 \\
		$f_{k,\gamma_5}$              & 0.029                    &  0.104                  &  0.034\\
		$f_{n,.\gamma_1}$            & 0.027   		            &  0.094  		    	  & 0.036 \\
		$f_{n,\gamma_2}$           & 0.026                    &  0.089			      & 0.032  \\
		$f_{n,hcvu}$               & 0.029                    &  0.085 		   		  & 0.035\\
		$f_{n,\gamma_4}$           & 0.024                    &  0.086                  & 0.026 \\
		$f_{n,\gamma_5}$            & 0.024                    &  0.089                  & 0.024 \\
		\hline      	
	\end{tabular}
\end{center}
	
\end{table}

\begin{table}[h]

	\caption{$L_2$ error over $100$ replicates using Epanechnikov's kernel for model (first column) and 5 (second column) in $\mathbb{R}^2$ with $\mu_1=(-1,1)$ and $\mu_2=(1,1)$. In both models, we used $2000$ uniformly distributed points for the test sample, in model 4 on $[-1,1]^2$ while in model 5 on $[-2,2]^2$. In both models, $k=l=2000$ and the measure of the $B(\epsilon,x)$ are estimated using $20000$ uniformly distributed points in $[-2,2]^2$.}
	\label{normal2}\tiny
	\begin{center}
		\begin{tabular}{|c|c|c|}	
			\hline
			&\multicolumn{2}{c|}{$\sigma_1^2=.5,\sigma_2^2=.3,\rho=.2$} \\                                
			\hline
			${\hat{f}_{\text{agg}}}$     & \textbf{0.071}           & \textbf{0.063}  \\ 
			$f_{k,\gamma_1}$              & 0.105                   &  0.113\\
			$f_{k,\gamma_2}$             & 0.103                    & 0.109 \\
			$f_{k,\gamma_3}$                  & 0.102                    & 0.106 \\
			$f_{k,\gamma_4}$             & 0.100                    & 0.103  \\
			$f_{k,\gamma_5}$              & 0.099                    & 0.101\\
			$f_{n,\gamma_1}$               & 0.096   		            & 0.094 \\
			$f_{n,\gamma_2}$              & 0.095                    &  0.092\\
			$f_{n,hcvu}$                 & 0.094                    & 0.107\\
			$f_{n,\gamma_4}$             & 0.093                    & 0.088  \\
			$f_{n,\gamma_5}$              & 0.093                    & 0.086\\
			\hline      	
		\end{tabular}
	\end{center}
	
\end{table}

The results in tables \ref{beta} to \ref{normal2} show that except for some results in Table \ref{beta}, the best performance is obtained by the aggregated estimator.  Moreover,
Table \ref{beta} also shows that in 5 over 8 models, this is also the case.

\begin{figure}[h]
	\begin{center}
		\includegraphics[scale=.03]{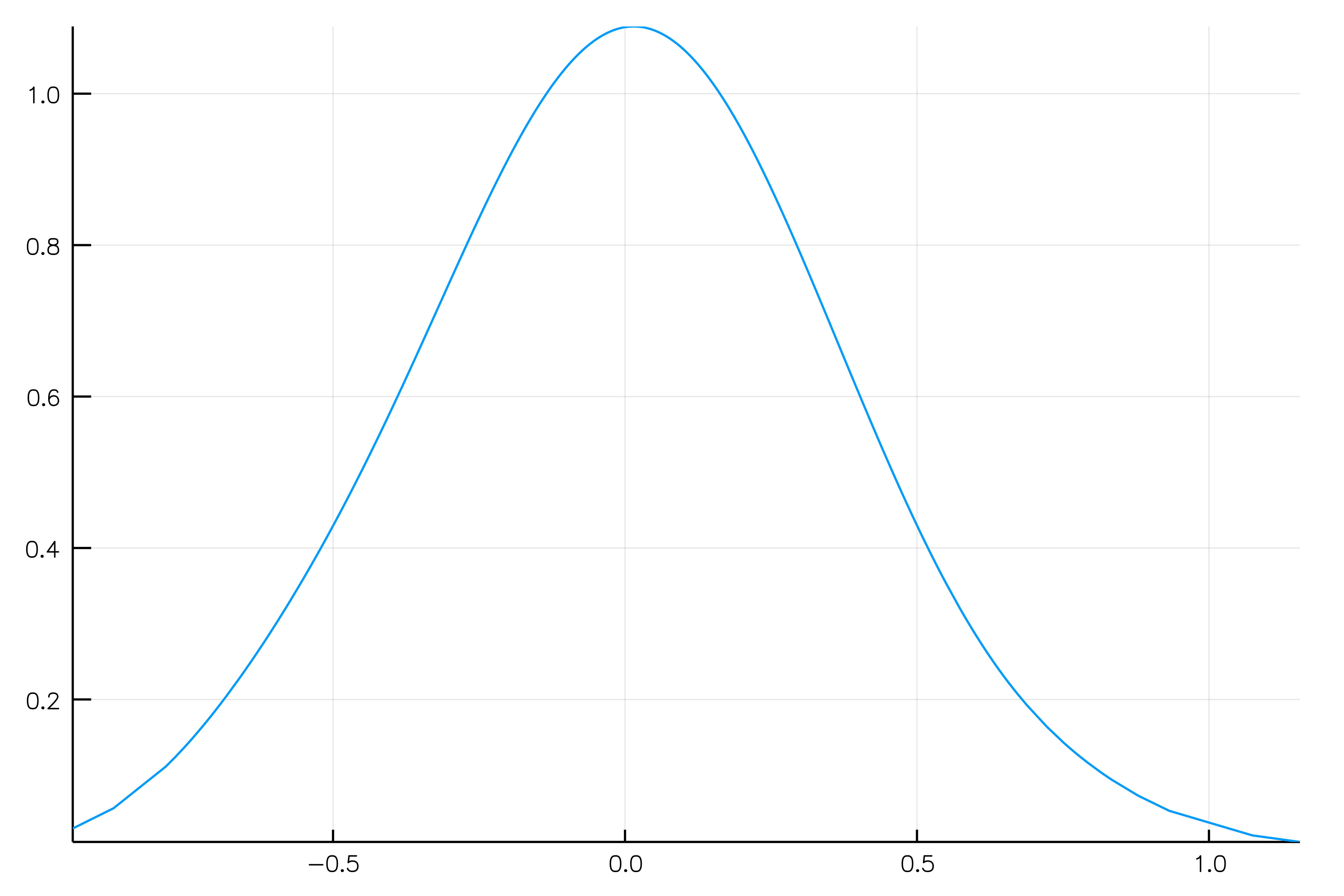} 
		\caption{Estimation of $N(0,f(x))$ for $x=(1/2,1/2)$ and $f$ the bi-variate normal.}
		\label{figth2}
	\end{center}
\end{figure}

To  illustrate Theorem \ref{tcl}, we have considered a bi-variate normal distribution with variance $\Sigma=Id$ and mean $(0,0)$. We fixed the point $x$ as $(1/2,1/2)$ for the normal distribution. The measure $\mu(B^*(\epsilon,x))$ is computed exactly from the density, in this case $f(1/2,1/2)=0.1239$ and for $\epsilon=0.005$, $\mu(B^*(\epsilon,x))=0.5088$. We have chosen $l=1000$, $k=6000$, computed $\sqrt{\mu(B^*(\epsilon,x))l}(\hat{f}_{\text{agg}}(x)-f(x))$ and repeated 1000 times. The estimator $\hat{f}_{\text{agg}}$ was built using $f_{k,hcv}$ (with Gaussian kernel). The density of the $N(0,f(x))$ was estimated using a kernel density estimator,  with a univariate Gaussian kernel with bandwidth 0.15. The result is shown in Figure \ref{figth2} and the summary is given in
Table \ref{summary}. The  p-value of Shapiro-Wilks test is $0.2772$ and $0.9978$ for Lilliefors test of normality.

\begin{table}[h]
\caption{Summary of the simulations for Theorem 2.}
\label{summary}
\begin{center}
\begin{tabular}{|c|c|c|c|c|c|c|}
\hline
Min.     & 1st. Qu  & Median  & Mean    & Var & 3rd Qu & Max\\
\hline
-0.946&-0.223 & 0.012& 0.010& 0.113 & 0.236&1.156 \\
\hline
\end{tabular}

\end{center}
\end{table}

\section{Final Remarks}
\begin{itemize}
	\item[1)] We have proposed a new non-linear aggregation method for density estimation and we have studied its asymptotic properties and limit distribution under quite mild assumptions.
	\item[2)]  The aggregated estimator behaves better than all of the density estimators used for the aggregation.
    \item[3)] We performed a small simulation study, which shows that in all cases the aggregation outperforms the kernel rules built with the sample $\mathcal{D}_k$. In addition, in most of them, it outperforms the kernel rules built with the whole sample $\mathcal{D}_n$. 
    \item[4)] Our simulations suggest that the second term in \eqref{prop1} is negligible with respect to the first term, but we were not able to prove this point theoretically. 
    \item[5)] The aggregation method is quite sensitive to the choice of the parameter $\epsilon$; however, because it is shown in the tables,  our recipe seems to work well.
\end{itemize}

\section{Appendix} \label{apA}

\textit{Proof of Proposition \ref{opt}}

We start by decomposing the objective function,
	\begin{multline*}
	\mathbb{E}|\hat{f}_{\text{agg}}(X)-f(X)|^2=\mathbb{E}|\hat{f}_{\text{agg}}(X)-T(\mathbf{f_k}(X))|^2
	+\mathbb{E}|T(\mathbf{f_k}(X))-f(X)|^2\\-2\mathbb{E}\Big[\big[\hat{f}_{\text{agg}}(X)-T(\mathbf{f_k}(X))\big]\big[T(\mathbf{f_k}(X))-f(X)\big]\Big].
	\end{multline*}
	Conditionally to $\mathbf{f_k}(X)$ and $\mathcal{D}_n$, $\hat{f}_{\text{agg}}(X)$ is constant, then
	\begin{multline*}\mathbb{E}\Big[\big[\hat{f}_{\text{agg}}(X)-T(\mathbf{f_k}(X))\big]\big[T(\mathbf{f_k}(X))-f(X)\big]\Big]=\\
	\mathbb{E}\Big[\hat{f}_{\text{agg}}(X)-T(\mathbf{f_k}(X))     \mathbb{E} \big[T(\mathbf{f_k}(X))-f(X)|\mathbf{f_k}(X),\mathcal{D}_n\big]\Big].
	\end{multline*}
	From $\sigma(\mathbf{f_k}(X),\mathcal{D}_n)\subset \sigma(\mathbf{f_k}(X))$ it follows 
	$ \mathbb{E}\big[\mathbb{E}[f(X)|\mathbf{f_k}(X)]|\mathbf{f_k}(X),\mathcal{D}_n\big]=$\break $\mathbb{E} \big[f(X)|\mathbf{f_k}(X),\mathcal{D}_n\big]$, then $ \mathbb{E} \big[T(\mathbf{f_k}(X))-f(X)|\mathbf{f_k}(X),\mathcal{D}_n\big]=0$,
	which implies that	$$\mathbb{E}\Big[\big[\hat{f}_{\text{agg}}(X)-T(\mathbf{f_k}(X))\big]\big[T(\mathbf{f_k}(X))-f(X)\big]\Big]=0.$$
	Lastly since $\mathbb{E}|T(\mathbf{f_k}(X))-f(X)|^2= \min_g \mathbb{E}[g(\mathbf{f_k}(X))-f(X)]^2$, where the minimum is taken over the functions $g$ such that $g(\mathbf{f}_k(X))\in L^2$, \eqref{prop1} follows.\\

\textit{Proof of Theorem \ref{thcons}}
First let us bound the second term in \eqref{prop1},	
		\begin{multline*}
		\mathbb{E}|\hat{f}_{\text{agg}}(X)-T(\mathbf{f_k}(X))|^2\leq \\
		2\mathbb{E}\Bigg[\frac{1}{l\mu(B(\epsilon,X))}\Big[\sum_{j=1}^l \mathbb{I}_{B(\epsilon,X)}(X_{k+j})-\mathbb{E}[\mathbb{I}_{B(\epsilon,X)}(X_{k+1})|\mathcal{D}_k,\mathbf{f_k}(X)]\Big]\Bigg]^2\\+
		2\mathbb{E}\Bigg[\frac{1}{\mu(B(\epsilon,X))}\mathbb{E}\Big[\mathbb{I}_{B(\epsilon,X)}(X_{k+1})|\mathcal{D}_k,\mathbf{f_k}(X)\Big]-T(\mathbf{f_k}(X))\Bigg]^2=I_1+I_2.
		\end{multline*}
		
		If we denote
		\begin{multline*}
		h(\mathbf{f_k}(X),\mathcal{D}_k)=\\ \mathbb{E}\Bigg[\Big[\sum_{j=1}^l \mathbb{I}_{B(\epsilon,X)}(X_{k+j})-l\mathbb{E}[\mathbb{I}_{B(\epsilon,X)}(X_{k+1})|\mathbf{f_k}(X),\mathcal{D}_k]  \Big]^2\Big|\mathbf{f_k}(X),\mathcal{D}_k\Bigg],
		\end{multline*}
		then 
		$$I_1=2\mathbb{E}\Bigg[\frac{1}{(l\mu(B(\epsilon,X)))^2}h(\mathbf{f_k}(X) \Bigg].$$
		
		Conditionally to $\mathbf{f_k}(X)$ and $\mathcal{D}_k$, the random variable $\sum_{j=1}^l \mathbb{I}_{B(\epsilon,X)}(X_{k+j})$ is binomial
		with probability $\mathbb{P}[{X_{k+1}\in B(\epsilon,X)}|\mathbf{f_k}(X),\mathcal{D}_k]$. Then
		\begin{align*}
		I_1=&\ 2\mathbb{E}\Bigg[\frac{1}{(l\mu(B(\epsilon,X)))^2}
		lP_X[{B(\epsilon,X)}|\mathbf{f_k}(X),\mathcal{D}_k](1-P_X[{B(\epsilon,X)}|\mathbf{f_k}(X),\mathcal{D}_k])\Bigg]\\
		\leq&\  2\mathbb{E}\Bigg[\frac{1}{(l\mu(B(\epsilon,X)))^2}lP_X[{B(\epsilon,X)}|\mathbf{f_k}(X),\mathcal{D}_k]\Bigg].\\
		\end{align*}
		We can bound $P_X[{B(\epsilon,X)}|\mathbf{f_k}(X),\mathcal{D}_k]\leq C\mu(B(\epsilon,X)),$ a.s., where $C=\sup f(x)$. Since $\{f_1\}_k,\dots,\{f_M\}_k$ are uniformly equicontinuous, for all $\epsilon>0$ there exists $\delta=\delta(\epsilon)>0$ such that for all $m=1,\ldots,M$, and all $k$, $|f_m(x)-f_m(y)|<\epsilon$ if $\|x-y\|<\delta$. By hypothesis we can assume that $\delta>\delta_0=\delta_0(\epsilon)$.  Then, $\mathcal{B}(x,\delta_0)\subset B(\epsilon,x)$, from where it follows,

		$$I_1\leq  2\mathbb{E}\Bigg[\frac{C}{l\mu(B(\epsilon,X))}\Bigg]\leq 
		2\frac{C}{\omega_d l\delta_0^d}\rightarrow 0\quad  \text{ as }l\rightarrow \infty.$$
		Regarding $I_2$, observe that,
 
		$$I_2=2\mathbb{E}\Bigg\{\frac{1}{\mu \left( B\left( \epsilon ,X\right) \right) }\int_{
			B\left( \epsilon ,X\right)} f(t)dt-
		T(\mathbf{f_k}(X))\Bigg\}^2. $$
	To prove this
		$\lim_{\epsilon \rightarrow 0}\lim_{l\rightarrow \infty }I_2 =0$, by Lemma \ref{L2conv},
		it is enough to show that 
		\[
		\lim_{\epsilon \rightarrow 0}\lim_{l\rightarrow \infty }\mathbb{E}\left( 
		\frac{1}{\mu\left( B\left( \epsilon ,X\right) \right) }\int_{
			B\left( \epsilon ,X\right) } f(t)dt-f(X)\right) ^{2}=0.
		\]%
		Because $f$ is bounded due to the dominated convergence theorem, it
		is enough to prove that 
		\begin{equation} \label{eq14}
		\lim_{\epsilon \rightarrow 0}\lim_{l\rightarrow \infty }\frac{1}{\mu
			\left( B\left( \epsilon ,X\right) \right) }\int_{B\left(
			\epsilon ,X\right)}f(t)dt=f(X) \quad a.s.
		\end{equation}
		
		Let $x$ such that for all $m=1,\dots,M$, $f_m(x)\rightarrow f(x)$. For such $x$ there exists $\epsilon_0=\epsilon_0(x)$ such that for all $\epsilon<\epsilon_0$, the sets $\overline{B(\epsilon,x)}$ are compact a.s, and $\overline{B^*(\epsilon,x)}$ is compact, then by Lemma \ref{lem} $\mu(B(\epsilon,x))\rightarrow \mu(B^*(\epsilon,x))$ a.s. By using again the dominated convergence theorem, we obtain that, with probability one,
		$$
		\lim_{\epsilon\rightarrow 0}\lim_{l\rightarrow \infty }\frac{1}{\mu \left( B\left( \epsilon ,X\right)
			\right) }\int_{B\left( \epsilon ,X\right)} f(t)dt-1=\\
		\lim_{\epsilon\rightarrow 0}\frac{1}{\mu \left( B^{\ast }\left( \epsilon ,X\right) \right) }%
		\int_{B^{\ast }\left( \epsilon ,X\right)}f(t)-f(X)dt.
		$$
		Lastly, \eqref{eq14} follows from the fact that for all $\epsilon>0$ and for all $t\in B^*(\epsilon,X)$, $|f(t)-f(X)|< \epsilon$.\\

\textit{Proof of Lemma \ref{L2conv}}\\

By Lemma 1.3 in \cite{alonso}, it is enough to prove that the sequence of $\sigma$-algebras $\{\sigma(f_i(X))\}_i$, $\mathbb{P}$-approaches $\sigma(f(X))$; i.e., for all $B\in \sigma(f(X))$ there exists $A_i\in \sigma(f_i(X))$ such that $\mathbb{P}(A_i\triangle B)\rightarrow 0$ as $i\rightarrow \infty$. Since $\sigma(f(X))=\sigma(\{f(X)^{-1}([a,b]),a,b\in\mathbb{R}\})$ is enough to consider $B=f(X)^{-1}([a,b])$ with $a<b$. Let us consider, for $\epsilon>0$, $B_i(\epsilon)=(f_i(X))^{-1}([a-\epsilon,b+\epsilon])$ and $H_i(\epsilon)=\cap^\infty_{j=i}\{\omega: |f_j(X(\omega))-f(X(\omega))|<\epsilon\}$.
	Let $A_i=\cap_{r=1}^\infty B_i(1/r)$, clearly $A_i\in \sigma(f_i(X))$. For all $\epsilon_1>0$,
	\begin{equation} \label{eq3}
	\mathbb{P}\big((\cap_{r=1}^\infty B_i(1/r))\triangle B\big)\leq \mathbb{P}\big(( \cap_{r=1}^\infty B_i(1/r))\triangle B\cap H_i(\epsilon_1)\big)+\mathbb{P}(H_i(\epsilon_1)^c).
	\end{equation}
	\begin{multline*}\mathbb{P}\big((\cap_{r=1}^\infty B_i(1/r))\triangle B\cap H_i(\epsilon_1)\big)=\mathbb{P}\big(( \cap_{r=1}^\infty B_i(1/r ))\cap B^c\cap H_i(\epsilon_1)\big)\ +\\
	\mathbb{P}\big((\cap_{r=1}^\infty B_i(1/r))^c\cap B\cap H_i(\epsilon_1)\big)=I_1+I_2.
	\end{multline*}
	Because the sequence of sets $\{B_i(1/r )^c\cap B\cap H_i(\epsilon_1)\}_r$ is increasing as $r$ increase,
	$$I_2=\mathbb{P}\Big(\cup_{r=1}^\infty \big(B_i(1/r )^c\cap B\cap H_i(\epsilon_1)\big)\Big)=\lim_{r\rightarrow \infty} \mathbb{P}\big(B_i(1/r)^c\cap B\cap H_i(\epsilon_1)\big),$$
	if $1/r <\epsilon_1$,
	\begin{align}\label{eq4}
	B_i(1/r )^c\cap B\cap H_i(\epsilon_1)\subset& \big\{\omega: f(X(\omega))\in [a,b]\setminus [a+\epsilon_1,b-\epsilon_1]\big\}.
	\end{align}
	
	Because the sequence of sets $\{B_i(1/r)\cap B^c\cap H_i(\epsilon_1)\}$ decreases as $r\rightarrow \infty$, 
	$$I_1=\mathbb{P}\big(\cap_{r=1}^\infty B_i(1/r)\cap B^c\cap H_i(\epsilon_1)\big)=\lim_{r\rightarrow \infty}\mathbb{P}\big(B_i(1/r)\cap B^c\cap H_i(\epsilon_1)\big),$$
	if $1/r<\epsilon_1$, $B_i(1/r)\cap H_i(\epsilon_1) \subset f(X)^{-1}([a-2\epsilon_1,b+2\epsilon_1]),$ 	then
	\begin{align*}
	\mathbb{P}\big(B_i(1/r)\cap B^c\cap H_i(\epsilon_1)\big)\leq &\  \mathbb{P}\big( B^c\cap   f(X)^{-1}([a-2\epsilon_1,b+2\epsilon_1])\big)\\
	=&\ \mathbb{P}(a-2\epsilon_1\leq f(X)<a)+\mathbb{P}(b< f(X)\leq b+2\epsilon_1).
	\end{align*}
	Because $\mathbb{P}(f(X)=a)=0=\mathbb{P}(f(X)=b)$, for all $\delta>0$ there exists $\epsilon_1>0$ such that 	$\mathbb{P}(a-2\epsilon_1\leq f(X)<a+\epsilon_1)+\mathbb{P}(b-\epsilon_1< f(X)\leq b+2\epsilon_1)<\delta.$
	By \eqref{eq4} and \eqref{eq3} for all $\delta>0$, $\mathbb{P}\big(   (\cap_{r=1}^\infty B_i(1/r)    )\triangle B\big)\leq \delta +\mathbb{P}(H_i(\epsilon_1)^c).$
	For all $\epsilon>0$, $\mathbb{P}(H_i(\epsilon)^c)\rightarrow 0$ as $i\rightarrow \infty$, from where it follows
	$\lim_{i\rightarrow \infty}\mathbb{P}\big((\cap_{r=1}^\infty B_i(1/r))\triangle B\big)=0.$

\textit{Proof of Lemma \ref{lem}}\\

Let us fix $x$ such that $\mu[B^*(\epsilon+\delta,x)\setminus B^*(\epsilon-\delta,x)]\rightarrow 0$ as $\delta\rightarrow 0$, and $f_m(x)\rightarrow f(x)$ for all $m$. First we will prove that, for all $\delta>0$, with probability one, for $k$ large enough $B(\epsilon,x)\subset B^*(\epsilon+\delta,x)$. Since $f$ and $f_m$ are uniformly continuous on $\overline{B(\epsilon,x)}$ we can take $\mathcal{Q}=\{q_1,\dots,q_s\}\subset B(\epsilon,x)$ such that for all $y\in B(\epsilon,x)$ there exists $q=q(y)\in\mathcal{Q}$ such that $|f(y)-f(q)|<\delta/3$.  Let $y\in B(\epsilon,x)$ and $q=q(y)\in \mathcal{Q}$ such that $|f(y)-f(q)|<\delta/3$. Then, 
	$$|f(y)-f(x)|\leq |f(y)-f(q)|+|f(q)-f(x)|.$$
	Meanwhile,
	$$|f(q)-f(x)|\leq |f(q)-f_m(q)|+|f_m(q)-f_m(x)|+|f_m(x)-f(x)|.$$
	Let $k$ large enough such that for all $q\in \mathcal{Q}$, $|f(q)-f_m(q)|<\delta/3$, and $|f_m(x)-f(x)|<\delta/3$. Then,  
	$|f(y)-f(x)|\leq \delta+\epsilon.$
	
	Now let us prove that for all $\delta>0$ such that $\epsilon-3\delta>0$,  $B^*(\epsilon-\delta,x)\subset B(\epsilon,x)$ a.s, as 
	$k\rightarrow \infty$. Proceeding as before, let us consider $\mathcal{Q'}=\{q'_1,\dots,q'_r\}\subset B^*(\epsilon-\delta,x)$ such that  for 
	all $y\in B^*(\epsilon-\delta,x)$ there exists $q'=q'(y)\in\mathcal{Q'}$ such that $|f_m(y)-f_m(q')|<\delta/3$ for all $m$. 
	Let $y\in B^*(\epsilon-\delta,x)$ and $q'\in \mathcal{Q'}$, such that $|f_m(x)-f_m(q')|<\delta/3$ for all $m$. Then,
	for all $m\in \{1,\dots,M\}$,
	$$|f_m(y)-f_m(x)|\leq |f_m(y)-f_m(q')|+|f_m(q')-f(q')|+|f(q')-f(x)|+|f(x)-f_m(x)|.$$
	
	Let $k$ be large enough such that for all $q'\in \mathcal{Q}'$ $|f(q')-f_m(q')|<\delta/3$, and $|f_m(x)-f(x)|<\delta/3$.
	Because $q'\in B^*(\epsilon-\delta,x)$, $|f(x)-f(q')|<\epsilon-\delta$, from where it follows that $y\in B(\epsilon,x)$. Lastly  
	$\mu[B^*(\epsilon+\delta,x)\setminus B^*(\epsilon-\delta,x)]\rightarrow 0$ implies \eqref{lem1}. To prove \eqref{lem2} let $\gamma>0$ and $\kappa$ small enough such that $P_X(B^*(\epsilon+\gamma,x))-P_X(B^*(\epsilon-\gamma,x))<\kappa$, for that $\delta$, with probability one, we can take $k$ large enough such that, $P_X(B^*(\epsilon-\gamma,x))\leq P_X(B(\epsilon,x))\leq P_X(B^*(\epsilon+\gamma,x)).$\\

\textit{Proof of Theorem \ref{th12}}

Let us denote $K_\epsilon(x)=K(x/\epsilon)$, then,
	\begin{multline*}
	\mathbb{E}|\tilde{f}_{\text{agg}}(X)-T(\mathbf{f_k}(X))|^2\leq\\
	2\mathbb{E}\Bigg\{\frac{1}{l\mu(B(\epsilon,X))}\Big[\sum_{j=1}^l 
	\prod_{m=1}^M K_\epsilon \big(f_m(X_{k+j})-f_m(X)\big)  
	\\ -l\mathbb{E}[ \prod_{m=1}^M K_\epsilon \big(f_m(X_{k+1})-f_m(X)\big)|
	\mathcal{D}_k,\mathbf{f_k}(X)]\Big]\Bigg\}^2\\+
	2\mathbb{E}\Bigg\{\frac{1}{\mu(B(\epsilon,X))}\mathbb{E}[ \prod_{m=1}^M
	K_\epsilon \big(f_m(X_{k+1})-f_m(X)\big)|\mathcal{D}_k,\mathbf{f_k}(X)]-
	T(\mathbf{f_k}(X))\Bigg\}^2\\
	=I_1+I_2.
	\end{multline*}
	Observe that 
	$$
	I_1=2\mathbb{E}\Bigg\{\frac{1}{(l\mu(B(\epsilon,X)))^2}\mathbb{V}\Big[\sum_{j=1}^l \prod_{m=1}^M K_\epsilon \big(f_m(X_{k+j})-f_m(X)\big)\Big|\mathbf{f_k}(X),\mathcal{D}_k\Big]\Bigg\}$$
	If we bound $K_\epsilon(x)\leq c_2$, then we get 
	$$I_1\leq 2\mathbb{E}\left[\frac{c_2^{2M}}{l\mu(B(\epsilon,X))^2}\right]. $$
	Proceeding as in Theorem \ref{thcons}, it is proved that $\lim_\epsilon\lim_l I_1=0$ a.s. 
	
	Regarding $I_2$ observe that,
 
	$$I_2=2\mathbb{E}\Bigg\{\frac{1}{\mu \left( B\left( \epsilon ,X\right) \right) }\int_{\mu \left(
B\left( \epsilon ,X\right) \right) }\prod_{m=1}^{M}K_{\epsilon }\left(
f_{m}(t)-f_{m}\left( X\right) \right) f(t)dt-
	T(\mathbf{f_k}(X))\Bigg\}^2. $$
  To prove that
	$\lim_{\epsilon \rightarrow 0}\lim_{l\rightarrow \infty }I_2 =0$, by Lemma \ref{L2conv},
	it is enough to show that 
\[
\lim_{\epsilon \rightarrow 0}\lim_{l\rightarrow \infty }\mathbb{E}\left[
\frac{1}{\mu \left( B\left( \epsilon ,X\right) \right) }\int_{\mu \left(
B\left( \epsilon ,X\right) \right) }\prod_{m=1}^{M}K_{\epsilon }\left(
f_{m}(t)-f_{m}\left( X\right) \right) f(t)dt-f(X)\right] ^{2}=0.
\]%
Because $f$ and $K$ \ are bounded, due to dominated convergence theorem, it
is enough to prove that 
\begin{equation} \label{eq11}
\lim_{\epsilon \rightarrow 0}\lim_{l\rightarrow \infty }\frac{1}{\mu
\left( B\left( \epsilon ,X\right) \right) }\int_{\mu \left( B\left(
\epsilon ,X\right) \right) }\prod_{m=1}^{M}K_{\epsilon }\left(
f_{m}(t)-f_{m}\left( X\right) \right) f(t)dt=f(X).
\end{equation}

Indeed, by using again dominated convergence theorem,  together with Lemma \ref{lem},   we obtain that,
\begin{multline*}
\lim_{\epsilon\rightarrow 0}\lim_{l\rightarrow \infty }\frac{1}{\mu \left( B\left( \epsilon ,X\right)
\right) }\int_{\mu \left( B\left( \epsilon ,X\right) \right)
}\prod_{m=1}^{M}K_{\epsilon }\left( f_{m}(t)-f_{m}\left( X\right) \right)
f(t)dt=\\
\lim_{\epsilon\rightarrow 0}\frac{1}{\mu \left( B^{\ast }\left( \epsilon ,X\right) \right) }%
\int_{\mu \left( B^{\ast }\left( \epsilon ,X\right) \right) }\left(
K_{\epsilon }\left( f(t)-f\left( X\right) \right) \right) ^{M}f(t)dt.
\end{multline*}
Meanwhile,

\begin{multline*}
\frac{1}{\mu \left( B^{\ast }\left( \epsilon ,X\right) \right) }\int_{\mu
\left( B^{\ast }\left( \epsilon ,X\right) \right) }\left( K_{\epsilon
}\left( f(t)-f\left( X\right) \right) \right) ^{M}f(t)dt-f(X)=
\\
\frac{1}{\mu \left( B^{\ast }\left( \epsilon ,X\right) \right) }\int_{\mu
\left( B^{\ast }\left( \epsilon ,X\right) \right) }\left( K_{\epsilon
}\left( f(t)-f\left( X\right) \right) \right) ^{M}\left( f(t)-f(X)\right)
dt+\\
-f(X)\Big[1-\frac{1}{\mu \left( B^{\ast }\left( \epsilon ,X\right) \right) }%
\int_{\mu \left( B^{\ast }\left( \epsilon ,X\right) \right) }\left(
K_{\epsilon }\left( f(t)-f\left( X\right) \right) \right) ^{M}dt\Big].
\end{multline*} 
 Lastly, \eqref{eq11} follows from \eqref{eq00} and the fact that for all $t\in B^*(\epsilon,X)$ $|f(t)-f(X)|\leq \epsilon$.\\

\textit{Proof of Theorem \ref{tcl}}

 First, let us prove that 
	\begin{equation} \label{eqconv}
	\lim_{l\rightarrow \infty}\lim_{k\rightarrow \infty} \frac{P_X(B(\epsilon,x))}{\mu(B(\epsilon,x))}=f(x) \quad \text{ a.s.}
	\end{equation}
	By Lemma \ref{lem} we get that for all fixed $l$,
	\begin{equation} \label{conv1}
	P_X(B(\epsilon,x))\rightarrow P_X(B^*(\epsilon,x))\quad  \text{ a.s.,}\quad \text{ as }k\rightarrow \infty
	\end{equation}
	and
	\begin{equation} \label{conv2}
	\mu(B(\epsilon,x))\rightarrow \mu(B^*(\epsilon,x))\quad  \text{ a.s.,}\quad \text{ as }k\rightarrow \infty.
	\end{equation}
	Lastly, from  
	\begin{equation} \label{eq5}
	\Big|\frac{P_X(B^*(\epsilon,x))}{\mu(B^*(\epsilon,x))}-f(x)\Big|\leq \frac{1}{\mu(B^*(\epsilon,x))}\int_{B^*(\epsilon,x)}|f(t)-f(x)|dt\leq \epsilon.
	\end{equation}
	it follows \eqref{eqconv}.  Let us write
	\begin{multline*}
	\sqrt{\mu(B^*(\epsilon,x))l }\bigg[\hat{f}_{\emph{agg}}(x)-f(x)\bigg]=\sqrt{\mu(B^*(\epsilon,x))l }\bigg[\hat{f}_{\emph{agg}}(x)- \frac{P_X(B(\epsilon,x))}{\mu(B(\epsilon,x))}\bigg]+\\
	\sqrt{\mu(B^*(\epsilon,x))l }\bigg[ \frac{P_X(B(\epsilon,x))}{\mu(B(\epsilon,x))}-f(x)\bigg]=I_1+I_2.
	\end{multline*}
	
	Since $f(x)>0$ then  $\mu(B^*(\epsilon,x))<\infty$ for $\epsilon$ small enough.
	From \eqref{conv1} and \eqref{conv2} together with  \eqref{eq5} and $l\epsilon^2\rightarrow 0$, it follows that $\lim_l\lim_k I_2=0$ a.s.\\

	Let us denote $Y_1,\dots,Y_l$ the random sample in $\mathcal{E}_l$ (i.e $X_{i+k}=Y_i$ for $i=1,\dots,l$). From \eqref{eq1} together with \eqref{conv1} and \eqref{conv2},
	\begin{equation*} 
	\lim_{k\rightarrow \infty} \ \ \hat{f}_{\text{agg}}(x)-\frac{P_X(B(\epsilon,x))}{\mu(B(\epsilon,x))}=\frac{\sum_{j=1}^l \mathbb{I}_{B^*(\epsilon,x)}(Y_j)}{l\mu(B^*(\epsilon,x))}-\frac{P_X(B^*(\epsilon,x))}{\mu(B^*(\epsilon,x))}\quad a.s.
	\end{equation*}
	Let us denote 
	$$f^*(x)=\frac{\sum_{j=1}^l \mathbb{I}_{B^*(\epsilon,x)}(Y_j)}{l\mu(B^*(\epsilon,x))}$$
	
	We will use the following version of the central limit theorem for triangular arrays.\\
	
	\textbf{Theorem} (\textbf{Lindeberg.}) 
	Let $Z_{l1},\dots,Z_{ll}$ independent r.v. such that for all $r=1,\dots,l$, $\mathbb{E}(Z_{lr})=m_{lr}$ and $\text{Var}(Z_{lr})=\sigma_{lr}^2<\infty$. Let us denote
	$V_l^2=\sum_{j=1}^l\sigma_{lj}^2.$
	If for all $\alpha>0$
	\begin{equation} \label{lindcond}
	\lim_{l\rightarrow+\infty}\frac{1}{V_l^2}\sum_{j=1}^l \mathbb{E}\left((Z_{lj}-m_{lj})^2\ind_{\{|Z_{lj}-m_{lj}|\geq \alpha V_l\}}\right)=0
	\end{equation}
	then
	\begin{equation}\label{lindth}
	\frac{1}{V_l}\sum_{j=1}^l (Z_{lj}-m_{lj})\stackrel{d}{\rightarrow} N(0,1)\quad \text{ as } l\rightarrow \infty.
	\end{equation}

	Let us consider $\epsilon=\epsilon_l\rightarrow 0$, define for $j=1,\dots,l$, 
	$$Z_{lj}=\frac{\mathbb{I}_{B^*(\epsilon,x)}(Y_{j})}{\mu(B^*(\epsilon,x))} \  \text{ then }  \ m_{lj}=\frac{P_X(B^*(\epsilon,x))}{\mu(B^*(\epsilon,x))}$$ 
	$$\text{ and  } \sigma_{lj}^2=\frac{P_X(B^*(\epsilon,x))(1-P_X(B^*(\epsilon,x)))}{\mu^2(B^*(\epsilon,x))},$$
	so 
	$$V_l=\sqrt{l\frac{P_X(B^*(\epsilon,x))(1-P_X(B^*(\epsilon,x)))}{\mu^2(B^*(\epsilon,x))}}.$$ 
	From $\mu\{y:f(x)=f(y)\}=0$ we get $P_X(B^*(\epsilon,x))\rightarrow 0$ as 
	$\epsilon\rightarrow 0$, and then using \eqref{eq5} it follows,  
	\begin{equation}\label{eq2}
	\lim_{l\rightarrow \infty} \frac{\sqrt{\mu(B^*(\epsilon,x))}}{\sqrt{P_X(B^*(\epsilon,x))(1-P_X(B^*(\epsilon,x)))}}=\sqrt{\frac{1}{f(x)}
	} \quad a.s..
	\end{equation}
	
	To prove \eqref{lindcond} 
	\begin{multline*}
	\frac{1}{V_l^2}\sum_{j=1}^l \mathbb{E}\left((Z_{l1}-m_{l1})^2\ind_{\{|Z_{lj}-m_{lj}|\geq \alpha V_l\}}\right)=\frac{l \mathbb{E}\Big((Z_{l1}-m_{l1})^2\ind_{\{|Z_{l1}-m_{l1}|\geq \alpha V_l\}}\Big)}{V_l^2}\\= \frac{1}{P_X(B^*(\epsilon,x))(1-P_X(B^*(\epsilon,x)))}\times \\  \mathbb{E}\left[\Big(\mathbb{I}_{B^*(\epsilon,x)}(Y_{1})-P_X(B^*(\epsilon,x))\Big)^2\ind_{|\mathbb{I}_{B^*(\epsilon,x)}(Y_{j})-P_X(B^*(\epsilon,x))|\geq \alpha \mu(B^*(\epsilon,x))V_l}\right].
	\end{multline*}
	Since $f$ is bounded and $\mu(B^*(\epsilon,x))l\rightarrow \infty $ it follows that $\mu(B^*(\epsilon,x))V_l\rightarrow \infty$. Then, with probability one, for $l$ large enough, $\ind_{|\mathbb{I}_{B^*(\epsilon,x)}(Y_{j})-P_X(B^*(\epsilon,x))|\geq \alpha \mu(B^*(\epsilon,x))V_l}=0$ and then it follows \eqref{lindcond}.

	Now from \eqref{lindth}, as $l\rightarrow \infty$,
	$$\frac{\mu(B^*(\epsilon,x))}{\sqrt{lP_X(B^*(\epsilon,x))(1-P_X(B^*(\epsilon,x)))}}\Bigg[\sum_{j=1}^l \Big[\frac{\mathbb{I}_{B^*(\epsilon,x)}(Y_j)}{\mu(B^*(\epsilon,x))}-\frac{P_X(B^*(\epsilon,x))}{\mu(B^*(\epsilon,x))}\Big]\Bigg]\stackrel{d}{\rightarrow} Z
	$$
	so, 
	$$\frac{\sqrt{\mu(B^*(\epsilon,x))}}{\sqrt{P_X(B^*(\epsilon,x))(1-P_X(B^*(\epsilon,x)))}} \sqrt{l\mu(B^*(\epsilon,x))}\Bigg[f^*(x)-\frac{P_X(B^*(\epsilon,x))}{\mu(B^*(\epsilon,x))}\Bigg]\stackrel{d}{\rightarrow} Z.$$
	where $Z=N(0,1)$. Then, from \eqref{eq2},
	$$\sqrt{l\mu(B^*(\epsilon,x))}\Bigg[f^*(x)-\frac{P_X(B^*(\epsilon,x))}{\mu(B^*(\epsilon,x))}\Bigg]\stackrel{d}{\rightarrow} N(0,f(x))\quad  \text{ as } l\rightarrow \infty.$$

\textit{Proof of Proposition \ref{prop2}}
	
		Let us denote $L(\lambda)=\{y:f(y)>\lambda\}$ the $\lambda$ level set of $f$, since $f$ is spherical for all $x\in int(S)$ being $S$ the support of $f$, $\nabla f(x)=2h'(\|x\|^2)x$ and then, using that $h'$ is a continuous function on a neighbourhood containing $\|x\|^2$, for $\epsilon$ small enough, $\nabla f(x)$ is a continuous function on $B^*(\epsilon,x)$. Since $f(x)>0$ $B^*(\epsilon,x)$ is bounded and then $f$ is Lipschitz on $B^*(\epsilon,x)$ for $\epsilon$ small enough. By Theorem 3.1 in \cite{fed59} we can write,
		\begin{equation} \label{area}
		\int_{B^*(\epsilon,x)}\|\nabla f(x)\|dx= \int_{-\epsilon}^\epsilon \mathcal{H}_{d-1}(\partial L(f(x)+t))dt,
		\end{equation}
		where $\mathcal{H}_{d-1}$ denotes the $(d-1)$-dimensional Hausdorff measure.
		Let us prove that $\mathcal{H}_{d-1}(\partial L(\lambda))$ is continuous for all $\lambda$ on a neighbourhood of $f(x)$. Observe that $\partial L(\gamma)=\{y:h(\|y\|^2)=\gamma\}$ implies $\mathcal{H}_{d-1}(\partial L(\gamma))=\mathcal{H}_{d-1}\big(\partial \mathcal{B}(0,\sqrt{h^{-1}(\gamma)})\big)$. Since $h$ is strictly decreasing there exists $h^{-1}$ (which is continuous on a neighbourhood of $\|x\|^2$ because $h$ is derivable) and $\|y\|^2=h^{-1}(\gamma)\rightarrow \|x\|^2$ as $\gamma\rightarrow f(x)$. By the Mean Value Theorem
		
		$$\int_{-\epsilon}^\epsilon \mathcal{H}_{d-1}(\partial L(f(x)+t))dt=2\epsilon \mathcal{H}_{d-1}(\partial L(\theta))\quad \text{ for some }f(x)-\epsilon\leq \theta\leq f(x)+\epsilon.$$
		Let us denote $M_\epsilon= \sup_{z\in B^*(\epsilon,x)}\|\nabla f(z)\|$ and $m_\epsilon= \inf_{z\in B^*(\epsilon,x)}\|\nabla f(z)\|$, then from \eqref{area}
		\begin{equation} \label{bound1}
		m_\epsilon \mu(B^*(\epsilon,x))\leq2\epsilon \mathcal{H}_{d-1}(\partial L(\theta)) \leq M_\epsilon \mu(B^*(\epsilon,x)).
		\end{equation}
		Since $h$ is decreasing we get that $B^*(\epsilon,x)=\{y:|h(\|y\|^2)-h(\|x\|^2)|<\epsilon\}$ decreases (to $\partial L(f(x))$) as $\epsilon$ decreases. From the continuity of $h'$ at $\|x\|^2$ it follows that $M_\epsilon= \sup_{z\in B^*(\epsilon,x)}\|\nabla f(z)\|=2\sup_{z\in B^*(\epsilon,x)}h'(\|z\|^2)\|z\|\rightarrow 2h'(\|x\|^2)\|x\|$ as $l\rightarrow \infty$. Analogously, $m_\epsilon\rightarrow 2h'(\|x\|^2)\|x\|$. Lastly, from the continuity of $\mathcal{H}_{d-1}(\partial L(\theta))$ and \eqref{bound1} we get that
		$$\lim_{l\rightarrow \infty} \frac{\mu(B^*(\epsilon,x))}{2\epsilon}= \frac{\mathcal{H}_{d-1}(\partial L(f(x)))}{\|\nabla f(x)\|}=\frac{2\pi^{d/2}\|x\|^{d-1}}{\Gamma(d/2)\|\nabla f(x)\|},$$
		where we have used that 
		$$\mathcal{H}_{d-1}(\partial \mathcal{B}(0,\|x\|))=\frac{2\pi^{d/2}\|x\|^{d-1}}{\Gamma(d/2)}.$$



\begin{thebibliography}{9}

\bibitem[\protect\citeauthoryear{Alonso and Brambila-Paz}{1998}]{alonso}
\rm{Alonso, A. and Brambila-Paz, F.} (1998).
$L^p$-Continuity of conditional expectations.
\textit{Journal of Mathematical Analysis and Applications} Vol. 221, pp. 161--176.

 \bibitem[\protect\citeauthoryear{Biau et al}{2016}]{biau:16}
\rm{Biau. G, Fischer, A. Guedj, B. and Malley, J.} (2016).
COBRA: A combined regression strategy.
\textit{Journal of Multivariate Analysis} Vol. 146, pp. 18-28.

\bibitem[\protect\citeauthoryear{Bellec}{2017}]{be:17}
\rm{Bellec, P. C.} (2017)
Optimal exponential bounds for aggregation of density estimators.
\textit{Bernoulli.} Vol 23(1), pp. 219--248.


\bibitem[\protect\citeauthoryear{Bourel and Ghattas}{2013}]{bg:13}
\rm{Bourel, M. and Ghattas, B.} (2013)
Aggregating density estimators: an empirical study.
\textit{Open Journal of Statistics}, Vol 3, pp. 334--355.

\bibitem[\protect\citeauthoryear{Breiman}{1996}]{brei:96}
\rm{Breiman, L.} (1996)
Bagging Predictors. \textit{Machine Learning} Vol. 24, No. 2, pp. 123-140.

\bibitem[\protect\citeauthoryear{Chac\'on and Duong}{2018}]{cd:18}
\rm{Chac\'on, J.E., and Duong, T.} (2018)
Multivariate Kernel Smoothing and Its Applications.
Chapman and Hall/CRC. ISBN 9781498763011.


\bibitem[\protect\citeauthoryear{Cholaquidis et al}{2016}]{ch:16}
\rm{Cholaquidis, A. Fraiman, R., Kalemkerian, J. and Llop, P.} (2016)
A nonlinear aggregation type classifier
\textit{Journal of Multivariate Analysis} Vol. 146, pp. 269--281.



\bibitem[\protect\citeauthoryear{Federer}{1959}]{fed59}
\rm{Federer, H.} (1959).
Curvature measures.
\textit{Trans. Amer. Math. Soc.}  Vol. 93, 418--491.

\bibitem[\protect\citeauthoryear{Fraiman et al}{1997}]{flm:97}
\rm{Fraiman, R., Liu, R. and Meloche, J.} (1997).
Multivariate density estimation by probing depth
\textit{$L_1$-Statistical Procedures and Related Topics}, IMS Lecture Notes -Monograph Series Vol. 31.



\bibitem[\protect\citeauthoryear{Lecu\'e}{2006}]{le:06}
\rm{Lecu\'e, G.} (2006)
Lower bounds and aggregation in density estimation.
\textit{Journal of Machine Learning Research.} Vol 7, pp. 971--981.

\bibitem[\protect\citeauthoryear{Einmahl and Mason}{2005}]{eima:05}
\rm{Einmahl, U. and Mason, D. M.} (2005)
Uniform in  bndwidth consistency of Kernel-type function estimators.
\textit{The Annals of Statistics}, Vol. 33, pp. 1380--1403.



\bibitem[\protect\citeauthoryear{Rigollet and Tsybakov}{2007}]{rt:07}
\rm{Rigollet, Ph. and Tsybakov, A.B.} (2007)
Linear and convex aggregation density estimators.
\textit{Mathematical Methods of Statistics} Vol. 16 (3) pp. 260–280.




\bibitem[\protect\citeauthoryear{Stute and Werner}{1991}]{stute}
\rm{Stute, W. and Werner, U.} (1991)
Nonparametric estimation of elliptically contoured densities.
\textit{Nonparametric functional estimation and related topics.}  Ed. G. Roussas, pp. 173-190. NATO ASI Series.




\end{thebibliography}


\end{document}